\documentclass[12pt]{article}
\newcommand\cF{\mathcal{F}}

\newcommand\cO{\mathcal{O}}

\newcommand\bC{\mathbf{C}}
\newcommand\bP{\mathbf{P}}

\newcommand\bZ{\mathbf{Z}}

\renewcommand\b{\beta}
\newcommand\ext{\mathrm{Ext}}

\newcommand\Hom{\mathrm{Hom}}
\newcommand\qed{\hfill QED}

\newcommand\spec{\mathrm{Spec}}

\newtheorem{prop}{Proposition}[section]
\newtheorem{conj}[prop]{Conjecture}
\newtheorem{defn}[prop]{Definition}
\newtheorem{lem}[prop]{Lemma}
\newtheorem{cor}[prop]{Corollary}

\begin{document}
  \begin{center}
Genus zero Gopakumar-Vafa invariants of contractible curves\\
Sheldon Katz\\
Departments of Mathematics and Physics\\
University of Illinois at Urbana-Champaign\\
Urbana, IL 61801 USA
  \end{center}

\bigskip\noindent
{\bf Abstract.} A version of the Donaldson-Thomas invariants of
a Calabi-Yau threefold is proposed as a conjectural mathematical definition
of the Gopakumar-Vafa invariants.  These invariants have a local
version, which is verified to satisfy the required properties
for contractible curves.  This provides a new viewpoint on the 
computation of the local Gromov-Witten invariants of contractible curves by
Bryan, Leung, and the author.

\section{Introduction.}
Let $X$ be a Calabi-Yau threefold, $\b\in H_2(X,\bZ)$, and $g$ a
nonnegative integer.  The Gopakumar-Vafa invariants $n^g_\b$ were
first introduced as an integer-valued index arising from D-branes and
M2-branes wrapping holomorphic curves in string theory and M-theory
\cite{gv}.  They are claimed to satisfy the remarkable identity
\begin{equation}
\label{gvgw}
\sum_{\b,g} N^g_\b
q^\b\lambda^{2g-2}=\sum_{\b,g,m}
\frac{n^g_\b}m\left(2\sin\left(\frac{m\lambda}2\right)\right)^{2g-2}q^{m\b},
\end{equation}
where the $N^g_\b$ are the Gromov-Witten invariants of $X$.  The coefficient
of $\lambda^{-2}$ in (\ref{gvgw}) yields the identity
\begin{equation}
\label{gvgw0}
N^0_\b=\sum_{k\mid\b}\frac{n_{\b/k}^0}{k^3}.
\end{equation}

There have been several attempts to provide a mathematical definition
including \cite{hst,kkv} but there is still no general direct
mathematical definition which passes all known tests. However, the
invariants can be defined recursively in terms of the Gromov-Witten
invariants via (\ref{gvgw}).  When defined in this way, it is only a
priori clear that the Gopakumar-Vafa invariants are rational.  The
{\em integrality conjecture\/} asserts that the Gopakumar-Vafa
invariants are integers.

In this note, we consider a variant of the Donaldson-Thomas invariants of
Calabi-Yau threefolds and conjecture that they satisfy the condition
(\ref{gvgw0}) required to be the genus~0
Gopakumar-Vafa invariants.  For this reason, we will call them the {\em
DT-GV invariants\/}.  As evidence, we show that a local version of these
invariants satisfies the desired property for local contractible
curves.

This note is organized as follows.  In Section~\ref{sec:definition} we
review the definition of the Donaldson-Thomas invariants and show that
they apply in the situation considered here
(Proposition~\ref{thomasapplies}).  Conjecture~\ref{gvconj} asserts
that these are the genus~0 Gopakumar-Vafa invariants.  In
Section~\ref{sec:contract}, we review the results of \cite{bkl} on the
local Gromov-Witten theory of contractible curves and their
implication for the associated Gopakumar-Vafa invariants.  We then
define the local DT-GV invariants and
compute them directly for contractible curves (Proposition~\ref{mainthm}).  
As a corollary, they
satisfy the genus~0 properties required of the Gopakumar-Vafa
invariants (Corollary~\ref{gvcont}).

\medskip
In this note, we understand the Calabi-Yau condition to mean that
$K_X$ is trivial.  We also refer to the holomorphic Casson invariants
of \cite{thomas} and their extensions as {\em Donaldson-Thomas
invariants\/}.  Note that we are using this terminology more broadly
than it was used in \cite{mnop1}, where it only referred to the
holomorphic Casson invariants of ideal sheaves.

\medskip
I'd like to thank J.\ Bryan for drawing my
attention to \cite{siebert} and for suggesting improvements to an
earlier draft.  This work was supported in part by NSF grant
DMS-02-44412.

\section{Definition of the Invariants.}
\label{sec:definition}

In this section, we discuss certain Donaldson-Thomas invariants of
Calabi-Yau threefolds $X$ and conjecture that they are the genus~0
Gopakumar-Vafa invariants of $X$.

The Donaldson-Thomas moduli space considered here was already studied
in \cite{hst}, where it was used to give a different proposal for the
Gopakumar-Vafa invariants of $X$.  It was conjectured that the genus~0
invariants introduced there coincide with the Donaldson-Thomas
invariants studied in this note.  Here, we are adopting a different
viewpoint by focusing on the genus~0 invariants, and treating the
Donaldson-Thomas invariants as the proposal for the Gopakumar-Vafa
invariants.

\medskip
We start by reviewing the construction of Simpson's moduli space of
semistable sheaves of pure dimension~1 \cite{simpson}.  Another useful
reference for this section is \cite{hl}.

\smallskip
Let $(X,L)$ be a polarized Calabi-Yau threefold, and let $\b\in H_2(X,\bZ)$.
We consider the Simpson moduli space of sheaves of $\cO_X$ modules
\begin{equation}
\label{simpsonmoduli}
S(X,\b)=\left\{F\ {\rm of\ pure\ dim\ }1\mid \mathrm{ch}_2(F)=\b,\ \chi(F)=1,\ F\ {\rm semistable\
wrt\ } L
\right\}.
\end{equation}

For the convenience of the reader, we recall the notion of stability
adapted to the present situation.  For a sheaf $F$ of dimension~1, its
{\em Hilbert polynomial\/} is
\[
P_F(n)= (L\cdot\mathrm{ch}_2(F))n+\chi(F),
\]
and its {\em reduced Hilbert
polynomial\/} is 
\begin{equation}
\label{redhilb}
p_F(n)=\frac{P_F(n)}{L\cdot\mathrm{ch}_2(F)}.
\end{equation}
A sheaf $F$ of pure dimension~1 is {\em (semi)stable\/} if for all nontrivial
proper subsheaves $G\subset F$ we have $p_G(n) < p_F(n)$ (resp.\ $p_G(n)\le
p_F(n)$ for semistable).  By (\ref{redhilb}) the semistability condition is
equivalent to
\begin{equation}
  \label{equivss}
\frac{\chi(G)}{L\cdot\mathrm{ch}_2(G)}
\le
\frac{\chi(F)}{L\cdot\mathrm{ch}_2(F)}.
\end{equation}
for all subsheaves $G$ of $F$.
For a subscheme $Z\subset X$ we usually write $P_Z(n)$ (resp.\ $p_Z(n)$)
in place of $P_{\cO_Z}(n)$ (resp.\ $p_{\cO_Z}(n)$).

Note that any $F\in S(X,\b)$ is necessarily stable.  To see this,
observe that for $G$ a subsheaf of $F$, (\ref{equivss}) reads
$\chi(G)/(L\cdot\mathrm{ch}_2(G))\le
1/(L\cdot\b)$.  But since $L\cdot\mathrm{ch}_2(G)\le L\cdot\b$, equality
can only hold if $\chi(G)=1$ and
$L\cdot\mathrm{ch}_2(G)= L\cdot\b$.  But in that case, we must
have $G=F$, and we conclude that $F$ is stable.  

Thus $S(X,\b)$ is a projective variety.

\begin{prop}
\label{thomasapplies}
There exists a
canonically determined perfect obstruction theory on $S(X,\b)$
of virtual dimension 0.
\end{prop}

\medskip\noindent {\em Proof.\/} The construction of the virtual
fundamental class in \cite[Theorem~3.30]{thomas} goes through with
little modification (i.e.\ the positive rank assumption made in
\cite{thomas} is not needed for the current application).  We only
have to show that $\dim\ext^3(F,F)$ is independent of $F\in S(X,\b)$.
But by Serre duality,
\[
\ext^3(F,F)\simeq\left(\Hom(F,F)
\right)^*.
\]
Stability implies simplicity, so that $\Hom(F,F)\simeq\bC$ is given by
the scalar multiplication maps.  This immediately
implies that $\dim\ext^3(F,F)=1$ for all $F\in S(X,\b)$, so 
we have a canonical perfect obstruction theory.

The virtual dimension is $\dim\ \ext^1(F,F)-\dim\ \ext^2(F,F)$, which
vanishes by Serre duality.  Alternatively, this is a symmetric
obstruction theory so that the expected dimension is 0
\cite{behrend}. \qed

\begin{defn}
  The {\em DT-GV invariants\/} are the Donaldson-Thomas
invariants $n_\b(X)=\deg[S(X,\b)]^\mathrm{vir}
\in\bZ$.
\end{defn}

The argument in \cite{thomas} applies without modification to show that 
the $n_\b(X)$ are deformation invariants of the polarized Calabi-Yau $X$.

\medskip
We conjecture that the $n_\b(X)$ are the genus~0 Gopakumar-Vafa invariants.
More precisely, let $N_\b(X)$ be the genus~0 Gromov-Witten invariants of $X$.

\begin{conj}
\label{gvconj}
\[
  N_\b(X)=\sum_{k\mid\b}\frac{n_{\b/k}(X)}{k^3}.
\]
\end{conj}

\noindent
{\bf Remark.}  We can similarly define a virtual fundamental class for
more general smooth threefolds $X$ whenever $\dim\mathrm{Ext}^3(F,F)$
is independent of $F\in S(X,\b)$.  For example, this condition holds
for any $\b$ and $L$ if $-K_X$ is very ample, and in this case it can
be computed by a simple modification of the computation in 
\cite[Lemma~1]{mnop2} that the virtual dimension is $D=\int_\b c_1(X)$.  Note
that this coincides with the virtual dimensions of Gromov-Witten
theory and Donaldson-Thomas theory \cite{mnop2}.  This leads to a
conjectural definition of the integer-valued invariants whose
existence was conjectured in \cite{pandharipandeicm}.  Details will
be given elsewhere.

\section{Contractible curves.}
\label{sec:contract}

In this section, we define a local version of the DT-GV invariants 
and compute them in the case of a small neighborhood of a contractible
$\bP^1$.  We then verify Conjecture~\ref{gvconj} in this situation.

We begin by recalling the geometry of contractible curves and their
local Gromov-Witten theory.

\subsection{Contractible curves and their local
GW invariants.}

Suppose that $C\simeq\bP^1$ is an analytically contractible curve in a
Calabi-Yau threefold $X$.  Then there is an analytic contraction map
$f:X\to Y$ such that $f(C)$ is a point $p\in Y$ which is a normal
singularity of $Y$, and $f$ induces an isomorphism between $X-C$ and
$Y-p$.  We can and will replace $X$ by a small neighborhood of $C$ and $Y$
be a small neighborhood of $p$.

By a lemma of Reid \cite{pagoda}, the generic hyperplane section
through $p$ is a surface $S$ with an isolated
rational double point at $p$, and the proper transform of $S$ is a partial
resolution ${\bar S}\to S$, i.e.\ $\bar{S}$ has at worst rational double
points, and the minimal resolution $\tilde{S}\to
S$ of $S$ factors through $\bar{S}\to S$. Note that $C\subset\bar{S}$.  It
follows that, shrinking $X$ and $Y$ if necessary,
there is a map
$\pi:X\to \Delta$, where $\Delta\subset\bC$ is a disc containing the
origin, such that $\bar{S} =\pi^{-1}(0)$.

The possible singularity types of $S$ and $\bar{S}$ were classified in
\cite{km}, and there are six cases.  The singularity $p\in S$ is
either an $A_1,\ D_4,\ E_6,\ E_7$, or $E_8$ singularity, and there are
two subcases of the $E_8$ case.  The curve $C\subset \bar{S}$ is
the exceptional curve of a partial resolution of $S$, associated to a
particular vertex of the corresponding Dynkin diagram which we call
the {\em distinguished vertex\/}. The surface $\bar{S}$ is
obtained from the minimal resolution $\tilde{S}$ of $p\in S$ by a map
$\psi:\tilde{S}\to\bar{S}$ which
blows down all of the exceptional curves of $\tilde{S}\to S$ except
for the one corresponding to the distinguished vertex of the Dynkin
diagram.  

The vertex can be uniquely specified by giving the
multiplicity of that vertex in the fundamental cycle of the minimal
resolution.  The cases listed as a pair (singularity type,
multiplicity) are
\[
(A_1,1),\ (D_4,2),\ (E_6,3),\ (E_7,4),\ (E_8,5),\ (E_8,6).
\]
More details are given in \cite{km} or the appendix to this note.

These situations can be distinguished by an invariant of Koll\'ar
called the {\em length\/} of $C$, which is the multiplicity of
$f^{-1}(m_p)$ at the generic point of $C$.  The length $\ell$ 
coincides with the multiplicity of the distinguished vertex in
each of these cases.

In \cite{bkl}, it was shown that the Gromov-Witten invariants of
$X$ have a well-defined contribution arising from stable maps to $X$ with
image contained in $C$, and this contribution was computed.

\medskip
To state the result of \cite{bkl}, let $\bar{S}$ be as above and let
$I=I_{C,\bar{S}}$ be the ideal
sheaf of $C$ in $\bar{S}$.  For each $i$ with
$1\le i\le \ell$, let $I^{(i)}$ be the saturation of $I^i$ 
i.e.\ the smallest ideal sheaf containing $I^i$ which
defines a subscheme of $X$ of pure dimension~1 (necessarily having
support $C$).  Let $C_i\subset X$ be the subscheme of $X$ defined by
$I^{(i)}$.  Then it was shown in \cite{bkl}
that $C_i$ is an isolated point of the
component of the Hilbert scheme of $X$ that it is contained in.
Letting $n_i$ be the multiplicity of this point, the result is

\begin{prop}{\cite[Theorem~1.5]{bkl}}
\label{bklprop}
  The contribution of $C$ to the genus~0
Gromov-Witten invariant $N_{d[C]}(X)$ is
\[
\sum_{k\mid d}\frac{n_{d/k}}{k^3}.
\]
\end{prop}

\noindent
{\bf Remark.} In \cite{bkl}, the local Gromov-Witten invariants were
computed for every genus and were shown to be completely determined by
these $n_d$ according to (\ref{gvgw}) with $n_d^0=n_d$ and $n^g_d=0$ for
$g>0$.

\bigskip
For later use, we give here some properties of the curves $C_i$.

\begin{lem}
\label{ideallem}
\ 

\begin{description}
\item{i.} The curve $C_i$ is the unique 1~dimensional subscheme of $\bar{S}$
supported on $C$ without embedded points, having multiplicity $i$ at
its generic point.  
\item{ii.} The curve $C_\ell$ coincides with the
scheme-theoretic inverse image of $p$ by $f$.
\item{iii.} For each $i$ with $1\le i\le\ell$ we have
$\chi(\cO_{C_{i}})=1$ and $H^1(\cO_{C_i})=0$.  
\item{iv.} The sheaf $\cO_{C_i}$ is stable for each $i$ with $1\le i\le \ell$.
\end{description}
\end{lem}

\noindent
{\em Proof:\/} Away from the singularities of $\bar{S}$, a scheme
supported on $C$ with no embedded points and multiplicity $i$
coincides with the scheme defined by $I^i$, so we are reduced to a
local question near the singularities of $\bar{S}$.  We let $A$ be the
local ring of $\bar{S}$ at a singular point, and we abuse notation
slightly by again denoting by $I$ the prime ideal of $A$ corresponding
to $C$. Letting $J$ be an ideal of $A$ with $J\subset I$ and $A/J$ of
pure dimension~1, we see that $J$ is $I$-primary. Then by localizing
at $I$ we conclude that if $J_I$ is an ideal of $A_I$ of multiplicity
$i$ then $J$ is the $i^{\scriptstyle\mathrm{th}}$ symbolic power of
$I$, which again by primary decomposition is the saturation of $I^i$.
This proves (i).

Note that the scheme $f^{-1}(p)$ is contained in $\bar{S}$.  So (ii)
will follow from (i) if we can show that $f^{-1}(p)$ has no
embedded points.  But a curve with embedded points has nonconstant
regular functions; pulling back via $\psi^*$ gives a nonconstant
regular function on the exceptional scheme of $\tilde{S}\to S$, a
contradiction.

Assertion (iii) is trivial for $i=1$, since $C_1=C\simeq\bP^1$.  For
$i>1$ we use the short exact sequences
\begin{equation}
\label{restmult}
  0\to I_{C_{j}}/I_{C_{j+1}}\to \cO_{C_{j+1}}\to \cO_{C_j}\to 0.
\end{equation}
Now we compute that $I_{C_{j}}/I_{C_{j+1}}\simeq\cO_C(-1)$ for $1\le
j\le \ell-1$.  To see this, first note that $I_{C_{j}}/I_{C_{j+1}}$ is a
torsion-free sheaf of rank~1 on $C$, hence locally free.  We compute
its degree by the method in \cite{morbirat}.  Let $\pi:\tilde{S}\to
\bar{S}$ be the contraction map used to produce the partial
resolution, and let $\tilde{C} \simeq C$ be the proper transform of
$C$ via $\psi$.  There is a map
\begin{equation}
  \label{pullback}
g:\psi^*(I_{C_{j}}/I_{C_{j+1}})\to I_{\tilde{C}}^j/I_{\tilde{C}}^{j+1}
\simeq \cO_C(2j).
\end{equation}
The map $g$ is an isomorphism away from the singularities of
$\bar{S}$.  So we can find the degree of $(I_{C_{j}}/I_{C_{j+1}})$ by
explicit computation of the order of vanishing of $g$ at the
singular points of $\bar{S}$.  Doing this explicitly in each case
using the local equations in \cite{morbirat}, we compute that the
degree is $-1$.  Some details of the computation are given in the
Appendix.

The first assertion of (iii) now follows from
(\ref{restmult}).  The second assertion follows since $p$ is a
rational singularity.

To prove assertion (iv), we let $Z$ be a nontrivial proper
subscheme of $C_i$. We
have to show that the reduced Hilbert polynomials satisfy $p_{I_{Z,C_i}}(n)
< p_{C_i}(n)$.

First we show this inequality for each $Z=C_j$ with $j<i$.  We have
\[
P_{I_{C_j,C_i}}(n)=P_{C_i}(n)-P_{C_j}(n)=(i-j)(L\cdot C)n
\]
by Lemma~{\ref{ideallem}.  
Then the needed reduced
Hilbert polynomial is
$p_{I_{C_j,C_i}}(n)=n$
so that 
\begin{equation}
\label{ineq}
p_{I_{C_j,C_i}}(n)<p_{C_i}(n)=n+\frac1{i(L\cdot C)}
\end{equation}
as required.

If $Z$ is 0~dimensional, then $\mathrm{ch}_2(I_Z)=\mathrm{ch}_2
(\cO_{C_i})$, so we just have to show that
$\chi_{I_Z}(n)<\chi_{C_i}(n)$.  This follows immediately from
$\chi_{I_Z}(n)=\chi_{C_i}(n)-\chi(\cO_Z)$.

If $Z$ has length $j$ at the generic point of $C$, then $Z\subset C_j$
and there is a short exact sequence
\[
0\to I_{Z,C_i}\to I_{C_j,C_i}\to \cO_Y\to 0
\]
for some zero-dimensional subscheme $Y\subset C$.  The desired inequality
\[
p_{I_{Z,C_i}}(n)<p_{C_i}(n)
\] 
follows 
from $\mathrm{ch}_2(I_Z)=\mathrm{ch}_2(C_j)$ and $\chi(I_{Z,C_i})=
\chi(I_{C_j,C_i})-\chi(\cO_Y)\le \chi(I_{C_j,C_i})$ which implies
$p_{I_{Z,C_i}}(n)\le p_{I_{C_j,C_i}}(n)$ together with (\ref{ineq}).
\qed

\bigskip\noindent
{\bf Remark.} An alternative proof of $\chi(\cO_{C_i})=1$ can be given
using the rationality of the singularity which implies that $h^1(\cO_{C_i})=0$
and the use of $\psi^*$ to show that $h^0(\cO_{C_i})=1$.

\subsection{Computation of the DT-GV invariants.}

Let $L$ be a polarization of $X$.  Let $F$ be a
stable sheaf of pure dimension~1 supported on $C$ with
$\mathrm{ch}_2(F)=d[C]$ and $\chi(F)=1$.  The rigidity lemma of
\cite{ckm} implies that every deformation of $F$ is supported on $C$
(not necessarily with the reduced structure).
Thus $S(X,d[C])$ has a connected component $S_C(X,d[C])$ consisting of sheaves
supported on $C$.  We can therefore define the contribution $n_d(X)$
of $C$ to our
genus~0 invariants as the degree of the part of $[S(X,d[C])]^\mathrm{vir}$
supported on $S_C(X,d[C])$.  

Comparing Proposition~\ref{bklprop} with Conjecture~\ref{gvconj}, it
is clear what needs to be proven, and in fact:

\begin{prop}
\label{mainthm}
  $n_d(X)=n_d$.
\end{prop}

Note in particular that the invariants $n_d(X)$ are independent of the
polarization chosen.

\bigskip\noindent {\em Proof.\/} We have an exact sequence
\begin{equation}
\label{restricts}
0\to \cO_X\to \cO_X\to \cO_{\bar{S}}\to 0
\end{equation}
since $I_{\bar{S},X}=\pi^*(I_{0,\Delta})$, which implies that
$I_{\bar{S},X}$ is trivial.

Let $F\in S_C(X,d[C])$.  Tensoring with (\ref{restricts}) we get an
exact sequence
\begin{equation}
\label{restftos}
F\to F\to F\otimes\cO_{\bar{S}}\to 0.
\end{equation}
By the stability of $F$, 
the first map in (\ref{restftos}) must be a scalar multiplication,
which is either 0 or an isomorphism.  But in the latter case, we would
conclude that $F\otimes\cO_{\bar{S}}=0$, which is impossible.  Thus
the first map is 0, and we conclude that $F\simeq
F\otimes\cO_{\bar{S}}$, i.e.\ that the scheme-theoretic support of $F$
is contained in $\bar{S}$.  Since the intersection of all
surfaces $S\subset Y$ containing $p$ is the reduced point $p$, it
follows that the support of $F$ is contained in the intersection of
all the possible surfaces $\bar{S}$, which is the scheme-theoretic
inverse image
$f^{-1}(p)$, i.e.\ the scheme $C_\ell$ by Lemma~\ref{ideallem}.

Next we show that $F$ must be isomorphic to one of the sheaves
$\cO_{C_i}$ for $1\le i\le\ell$.  The argument is similar to the
argument in \cite{hst}.  By Lemma~\ref{ideallem}, we already know
that $\cO_{C_i}\in S(X,i[C])$ for $1\le i\le\ell$.

Since $\chi(\cO_F)=1$ and $h^2(\cO_F)=0$, it follows that $F$ has a
section, i.e.\ we have a map $s:\cO_X\to F$.  Since $F$ is an
$\cO_{C_\ell}$-module, the kernel of $s$ is the ideal sheaf of a
subscheme $Z\subset C_\ell$.  Since $F$ has pure dimension~1, it
follows that $Z$ has pure dimension~1 as well.  Lemma~\ref{ideallem}
implies that $Z\simeq C_i$ for some $i$ between 1 and $\ell$ and we
have an injection $\cO_{C_i}\hookrightarrow F$.

The Hilbert polynomial of $F$ is $P_F(n)=d(L\cdot[C])n+1$,
and the Hilbert polynomial of $C_i$ is $P_{C_i}(n)=
i(L\cdot[C])n+1$.  By stability, we conclude that $d\le i$, while $i\le d$
since $\cO_{C_i}$ is a subsheaf of $F$.  It follows that $i=d$, whence
$\cO_{C_i}\hookrightarrow F$ is an isomorphism.

We have shown that $S(X,d[C])$ consists of the single sheaf
$\cO_{C_d}$ if $d\le \ell$ and is zero if $d>\ell$, so we may as well
assume that $d\le\ell$.  Next we identify $S(X,d[C])$ with a component
of the Hilbert scheme.

Let $A$ be a local ring and let $\cF$ be a coherent sheaf on
$X\times \spec(A)$, flat over $\spec(A)$. Let $\rho:X\times
\spec(A)\to \spec(A)$ be the projection.  We have $H^1(\cO_{C_d})=0$
and $h^0(\cO_{C_d})=1$ by Lemma~\ref{ideallem}, so by standard base
change results, $\rho_*(\cF)$ is invertible, hence free of rank~1 on
$\spec(A)$.  The adjoint map to any isomorphism
$\cO_{\spec(A)}\to\rho_*(\cF)$ gives a map
\begin{equation}
  \label{surjection}
\cO_{X\times\spec(A)}\to\cF,
\end{equation}
necessarily a surjection, so that $\cF$ can be identified with a flat family 
of subschemes of $X$.  The remaining details are straightforward and left
to the reader.

The degree of the virtual fundamental class can be computed by
\cite[Theorem~4.6]{siebert}.  Since $S_C(X,d[C])$ is a point, the conclusion
is that it coincides with the degree of the Fulton Chern class, which
in turn is equal to the multiplicity of $C_d$ as a point
of $S_C(X,d[C])$ by \cite{fulton}. Alternatively, the methods of
\cite{bf} apply.  \qed

\begin{cor}
  \label{gvcont}
The local version of Conjecture~\ref{gvconj} holds for contractible
curves.
\end{cor}

\bigskip\noindent {\bf Remarks.}  

\smallskip\noindent (i) There is an alternative simpler proof of
Proposition~\ref{mainthm}.  The cited result of \cite{bkl} was proven
by showing that in a generic deformation of $X$, each $C_i$ deforms to
$n_i$ isolated curves, each isomorphic to $\bP^1$ with normal bundle
$\cO(-1)\oplus\cO(-1)$, then using the deformation invariance of the
Gromov-Witten invariants.  We can also invoke a local version of the
deformation invariance of the DT-GV invariants in combination with the
computation of the moduli space of stable sheaves supported on a
$\bP^1$ with normal bundle $\cO(-1)\oplus\cO(-1)$ done in \cite{hst}.
We have not gone this route to highlight the direct computability of the
DT-GV invariants.

\smallskip\noindent (ii) Since the genus $g$ Gopakumar-Vafa invariants
can be roughly thought of as the ``virtual number of genus $g$
Jacobians'' contained in the moduli space $S(X,\b)$, any reasonable
mathematical definition of the Gopakumar-Vafa invariants $n^g_\b$
should produce zero for $g>0$ whenever $S(X,\b)$ is zero-dimensional.
In particular, we should get zero for the higher genus local
Gopakumar-Vafa invariants of contractible curves, as computed by
Gromov-Witten theory in \cite{bkl}.  However, we will not single out a
specific proposed definition here.

\medskip\noindent
{\em Question:\/} For a general Calabi-Yau threefold $X$, are the DT-GV
invariants independent of $L$? 

\section*{Appendix:  Some computations.}
\label{computation}

In this appendix, we give some of the calculations supporting
the proof of part Lemma~\ref{ideallem}, (iii). The
illustrative examples given here should suffice to allow the interested
reader to carry the calculation to its completion.

First, we list the singularities of $\bar{S}$, which are visible from
Figure~1 in \cite{km}.

\[
  \begin{array}{|c|c|c|}\hline
    p&\ell&\mathrm{Sing}(\bar{S})\\ \hline
A_1&1&{\rm none}\\ \hline
D_4&2&A_1,\ A_1,\ A_1\\ \hline
E_6&3&A_2,\ A_2,\ A_1\\ \hline
E_7&4&A_3,\ A_2,\ A_1\\ \hline
E_8&5&A_4,\ A_3\\ \hline
E_8&6&A_4,\ A_2,\ A_1\\ \hline
  \end{array}
\]

We need to compute the order of vanishing of (\ref{pullback}) at the
singularities of $\bar{S}$.

A uniform treatment can be given for all cases except the $A_4$
singularity in the length~5 case.  Suppose that we are in any of the
other situations, with an $A_k$ singularity.  Then by the computations
in the appendix to \cite{morbirat}, we can describe $\bar{S}$ locally
as the hypersurface $xy+z^{k+1}=0$ with $C\subset\bar{S}$ having ideal
$(x,z)$.  We can then compute each $I^{(j)}$ locally and from that the
order of vanishing of the map $g$ of (\ref{pullback}) by looking at a
generator of $I^{(j)}/I^{(j+1)}$.

We illustrate with the $A_1$ case.  The blowup map $\psi:\tilde{S}\to
\bar{S}$ can be described on an affine piece of $\tilde{S}$ by
$(u,v)\mapsto(u^2v,v,uv)$.  Here $u=0$ defines the proper transform 
$\tilde{C}$ of $C$
and $(u,v)=(0,0)$ is the point of $\tilde{C}$ lying over the singularity
of $\bar{S}$, so the order of vanishing of a function is just given by
the exponent of $v$.  We get
\[
\begin{array}{|c|c|c|c|}\hline
  j&I^{(j)}&{\rm generator}&\mathrm{ord}(g)\\ \hline
1 & (x,z)&z&1\\ \hline
2 & (x)&x&1\\ \hline
3 & (x^2,xz)&xz&2\\ \hline
4 & (x^2)&x^2&2 \\ \hline
5 & (x^3,x^2z)&x^3&3\\ \hline
\end{array}
\]

So if for example $p$ is a $D_4$ singularity, we have three $A_1$ 
singularities.  Putting $j=1$ in (\ref{pullback}) we see that $g$
vanishes simply at each singularity, hence the degree of $I/I^{(2)}$ is
$2-3(1)=-1$ as claimed.

In the $A_2$ case, we similarly compute that $g$ has vanishing orders
$1,2,2,3,4$ in the respective cases $j=1,2,3,4,5$.

So if $p$ is an $E_6$ singularity, we have two $A_2$ singularities
and an $A_1$ singularity.  For $j=1$ we get
\[
\mathrm{deg}\left(\frac{I}{I^{(2)}}\right)=2-2(1)-1=-1
\]
while for $j=2$ we get
\[
\mathrm{deg}\left(\frac{I^{(2)}}{I^{(3)}}\right)=4-2(2)-1=-1
\]
as claimed.

\smallskip
The other cases are similar.

\end{document}